\documentclass{amsart}
\usepackage{graphicx}
\newtheorem{theorem}{Theorem}

\newtheorem{definition}[theorem]{Definition}
\newtheorem{example}[theorem]{Example}
\newtheorem{remark}[theorem]{Remark}

\def\QED{\quad\blackslug\lower 8.5pt\null}

\begin{document}

\vspace*{30mm}

\begin{center}
{\large \bf Geometry of Projective Planes
over Two-Dimensional Algebras}
\end{center}

\vspace*{3mm}

\begin{center}
{\large
 Maks A. Akivis and  Vladislav V. Goldberg}
\end{center}

\vspace*{3mm}

{\footnotesize
\noindent
\textit{$2000$ Mathematics Subject Classification}.
 Primary 53A20, Secondary 14M99.
} \newline

{\footnotesize
\textbf{Abstract.}
The authors study smooth lines on projective planes over
the algebra $\mathbb{C}$ of complex numbers, the algebra
$\mathbb{C}^1$ of double numbers, and the algebra $\mathbb{C}^0$
of dual numbers. In the space $\mathbb{R} P^5$, to these smooth lines
there correspond families of straight lines describing
point three-dimensional tangentially degenerate submanifolds $X^3$
of rank 2. The authors study focal properties of these submanifolds
and prove that they represent examples of different types of
tangentially degenerate submanifolds. Namely, the submanifold $X^3$,
corresponding in  $\mathbb{R} P^5$ to a smooth line $\gamma$ of
the projective plane $\mathbb{C}P^2$, does not have real singular points,
the submanifold $X^3$, corresponding in  $\mathbb{R} P^5$ to a smooth
line $\gamma$ of  the projective plane $\mathbb{C}^1 P^2$, bears two plane
singular lines, and finally the submanifold $X^3$, corresponding in
 $\mathbb{R} P^5$ to a smooth line $\gamma$ of  the projective plane
 $\mathbb{C}^0 P^2$, bears one singular line.
}

\setcounter{equation}{0}

\vspace*{5mm}

\textbf{1. Introduction.}
The theory of projective planes over algebras is the subject
belonging to the geometry and the algebra, and this subject
attracts the attention of both, algebraists and geometers.
This theory was considered in Pickert's book [Pi 75], and
in the separate chapters of the books  [B 70], [R 66], and
[R 97].

However, not so much is known about
the differential geometry of such projective planes.
Some questions in this direction were considered
in the paper [A 87]. In this paper the author studied
smooth lines in projective planes over the matrix algebra
and over some of its subalgebras. In this study
he used the mapping of the projective plane $M P^2$ over
the algebra $M$ of $(n \times n)$-matrices onto the
Grassmannian $G (n-1, 3n-1)$ of subspaces of dimension
$n-1$ of a real projective space $\mathbb{R} P^{3n-1}$.

In the current paper, we continue the investigations of [A 87].
However, we restrict ourselves to the study of smooth lines
on projective planes over the algebra $\mathbb{C}$ of complex
numbers, the algebra $\mathbb{C}^1$ of double numbers,
and the algebra $\mathbb{C}^0$ of dual numbers
 (for description of these algebras,
 see, for example, [Pa 63] or [Sc 66] or [R 97]).
In the space $\mathbb{R} P^5$, to these smooth lines
there correspond families of straight lines describing
point three-dimensional tangentially degenerate submanifolds $X^3$
of rank 2 (see [AG 93], Ch. 4).
We study focal properties of these submanifolds
and prove that they represent examples
of different types of tangentially degenerate
submanifolds. Namely, the submanifold $X^3$, corresponding in
 $\mathbb{R} P^5$ to a smooth line $\gamma$ of
 the projective plane $\mathbb{C}P^2$, does not have
real  singular points, the submanifold $X^3$, corresponding in
 $\mathbb{R} P^5$ to a smooth line $\gamma$ of
 the projective plane $\mathbb{C}^1 P^2$, bears two plane
   singular lines, and finally the submanifold $X^3$, corresponding in
 $\mathbb{R} P^5$ to a smooth line $\gamma$ of
 the projective plane $\mathbb{C}^0 P^2$, bears one double singular line.

In the last case, the submanifold $X^3$ is a generalization
of  tangentially degenerate hypersurfaces without singularities
in the four-dimensional Euclidean space $R^4$ constructed by
 Sacksteder in [Sa 60] and recently  considered by Bourgain.
 Note that Bourgain's hypersurface was described in
 the papers [W 95] and [I 98, 99a, 99b], and the authors
 proved (see [AG 00a]) that the hypersurfaces of  Sacksteder
and Bourgain are locally equivalent.

  Note also that it follows from the paper [A 87] that
  in the  projective plane $M P^2$ over the algebra
  of $(2 \times 2)$-matrices, there is no smooth lines
  different from straight lines.
  A family of straight lines in $\mathbb{R} P^5$ corresponding
  to those straight lines is the Grassmannian $G (1, 3)$
  of straight lines lying in a three-dimensional subspace
  of the space $\mathbb{R} P^5$.

Note that in [AG 00b], the authors found the basic types of
tangentially degenerate submanifolds and proved the structure
theorem for such submanifolds: an arbitrary  tangentially
degenerate submanifold is either irreducible or if it is
reducible, it is foliated into submanifolds of basic types.  The
finding of examples  of tangentially degenerate submanifolds of
basic and not basic types is important. Such examples can be found
in [AG 93, 00a, b], [AGL], [GH 79], [I 98, 99a, b], [L 99], [W
95]. In particular,  in [I 99a],
Ishikawa found real algebraic cubic nonsingular
tangentially degenerate hypersurface in $\mathbb{R} P^n$ for
$n = 4, 7, 13, 25$. These hypersurfaces have the structure of
homogeneous spaces of groups $\mathbf{SO} (3),
\mathbf{SU} (3), \mathbf{Sp} (3)$, and $F_4$,
respectively, and their projective  duals are linear
projections of  Veronese embeddings of projective planes
$\mathbf{K} P_2$ for $\mathbf{K} = \mathbf{R, C, H, O}$,
where $\mathbf{O}$ is the Cayley's octonions.

The examples we have constructed in this paper are
of the same nature as Ishikawa's examples but they are much
simpler.

\textbf{2. Two-dimensional algebras and their representation.}
There are known three two-dimensional algebras: the algebra
of complex numbers $z = x + i y$, where $i^2 = - 1$, the algebra
of double (or split complex) numbers $z = x + e y$, where $e^2 = 1$, and the
 algebra of dual numbers $z = x + \varepsilon y$, where
 $\varepsilon^2 = 0$. Here everywhere $x, y \in \mathbb{R}$.
 Usually these three algebras are denoted by $\mathbb{C},
 \mathbb{C}^1$, and $\mathbb{C}^0$ (see [R 97],  \S 1.1).
 These algebras are commutative and associative, and
 any  two-dimensional algebra is isomorphic to one
 of them.

 Each of these three algebras admits a  representation by means
 of the $(2 \times 2)$-matrices:
\begin{equation}\label{1}
z = x + i y \rightarrow \left(
\begin{array}{rr}
x & -y \\
y & x \\
\end{array} \right),
\end{equation}
\begin{equation}\label{2}
z = x + e y \rightarrow \left(
\begin{array}{ll}
x & y \\
y & x \\
\end{array} \right),
\end{equation}
and
\begin{equation}\label{3}
z = x + \varepsilon y \rightarrow \left(
\begin{array}{ll}
x & 0\\
y & x \\
\end{array} \right).
\end{equation}
In what follows, we will identify
the algebras  $\mathbb{C},  \mathbb{C}^1$, and $\mathbb{C}^0$
with their matrix representations.

The algebras  $\mathbb{C},  \mathbb{C}^1$, and $\mathbb{C}^0$
are subalgebras of the total matrix algebra $M$
formed by all  $(2 \times 2)$-matrices
\begin{equation}\label{4}
\renewcommand{\arraystretch}{1.3}
\left(
\begin{array}{ll}
x_0^0 & x_1^0 \\
x_0^1 & x_1^1
\end{array} \right),
\renewcommand{\arraystretch}{1}
\end{equation}
which is associative but not commutative.

The algebra  $\mathbb{C}$ does not have
zero divisors while the algebras
$\mathbb{C}^1, \mathbb{C}^0$, and $M$ have such divisors.
In the matrix representation, zero divisors of these algebras are
determined by the condition
$$
\renewcommand{\arraystretch}{1.3}
\det \left(
\begin{array}{ll}
x_0^0 & x_1^0 \\
x_0^1 & x_1^1
\end{array} \right) = 0.
\renewcommand{\arraystretch}{1}
$$
For the algebra $\mathbb{C}^1$ the last condition takes
the form
$$
x^2 - y^2 = 0,
$$
for the algebra $\mathbb{C}^0$
the form $x = 0$, and for the algebra $M$ the form
\begin{equation}\label{5}
x_0^0 x_1^1 - x_0^1 x_1^0 = 0.
\end{equation}

The elements of the algebras $\mathbb{C}^1$ and
$\mathbb{C}^0$, as well as the regular complex numbers (the elements
of the algebra $\mathbb{C}$), can be represented by the points
on the plane $x O y$. In this representation, the zero divisors
of the algebra $\mathbb{C}^1$ are represented by the points of the
straight lines $y = \pm x$, and the zero divisors
of the algebra $\mathbb{C}^0$ by the points of the $y$-axis.

The elements of the algebra $M$ are represented by the points
of a four-dimensional vector space, and its zero divisors by the
points of the cone (5) whose signature is $(2,2)$. Thus, to
 the algebra $M$ there corresponds a four-dimensional
 pseudo-Euclidean space $R_2^4$ of signature 2 with the isotropic
 cone (5). This cone bears two family of plane generators defined
 by the equations
\begin{equation}\label{6}
\frac{x_0^0}{x_0^1} = \frac{x_1^0}{x_1^1} = \lambda, \;\;
\frac{x_0^0}{x_1^0} = \frac{x_0^1}{x_1^1} = \mu,
\end{equation}
where $\lambda$ and $\mu$ are real numbers.

\textbf{3. The projective planes over the algebras
$\mathbb{C}, \mathbb{C}^1$,  $\mathbb{C}^0$, and
$\boldsymbol{M}$.}
Denote by $A$ one of the  algebras
$\mathbb{C},  \mathbb{C}^1$, and $\mathbb{C}^0$, or
$\boldsymbol{M}$ and consider a projective plane $AP^2$
over the algebra $A$ (see [B 70]). A point $X \in AP^2$
has three  matrix coordinates $X^0, X^1, X^2$ that
have respectively the form (1), (2), (3), or (4).
Since it is convenient to write point coordinates
 as a column-matrix, we write
\begin{equation}\label{7}
X = (X^0, X^1, X^2)^T.
\end{equation}
The matrix $X$ in (7) has 6 rows and 2 columns.
Of course, the columns of this matrix must be
linearly independent. The coordinates
$X^\alpha, \, \alpha = 0, 1, 2,$ are
defined up to a multiplication from the right
by an element  $P$ of the algebra $A$
which is not a zero divisor.

The columns of the matrix $X$
can be considered as coordinates of the points
$x_0$ and $x_1$ of a projective space  $\mathbb{R} P^5$,
and to the matrix $X$ there corresponds
the straight line $x_0 \wedge x_1$
in the space  $\mathbb{R} P^5$. So, we can
set $X = x_0 \wedge x_1$. The set of all straight lines of
the space  $\mathbb{R} P^5$ forms the Grassmannian $G (1, 5)$,
whose dimension, $\dim \, G (1, 5) = 2 \cdot 4 = 8$.
We express this in more detail. For all three planes
 $\mathbb{C} P^2,  \mathbb{C}^1 P^2$, and $\mathbb{C}^0 P^2$,
the coordinates of the point $x_0$ are the same:
$$
x_0 = (x^0_0, x^1_0, x^2_0, x^3_0, x^4_0, x^5_0)^T,
$$
and the coordinates of the point $x_1$ have different forms.
For the plane   $\mathbb{C} P^2$, they are
$$
x_1 = (-x^1_0, x^0_0, -x^3_0, x^2_0, -x^5_0, x^4_0)^T,
$$
for the plane $\mathbb{C}^1 P^2$, they are
$$
x_1 = (x^1_0, x^0_0, x^3_0, x^2_0, x^5_0, x^4_0)^T,
$$
and for the plane $\mathbb{C}^0 P^2$, they are
$$
x_1 = (0, x^0_0, 0, x^2_0, 0, x^4_0)^T.
$$

A straight line in the  plane $AP^2$
is defined by the equation
$$
U_0 X^0 + U_1 X^1 + U_2 X^2 = 0,
$$
where $U_\alpha, \, \alpha = 0, 1, 2,$ are matrices belonging to the
 algebra $A$.

Two points $X$ and $Y$ are called \textit{adjacent} if more that
one straight line passes through them. If
$$
 X = (X^0, X^1, X^2)^T, \;\;
Y =  (Y^0, Y^1, Y^2)^T
$$
are adjacent points, then the rank of the $(6 \times 4)$-matrix composed
of the matrix coordinates of $X$ and $Y$ is less than 4. If the rank
of this matrix is 4, then through the points $X$ and $Y$
there passes a unique straight line.

Under the mapping of  the plane $P^2 (M)$
into the space  $\mathbb{R} P^5$, to adjacent points there
correspond \textit{intersecting} straight lines.
If $X$ and $Y$ are not adjacent points, then in
 $\mathbb{R} P^5$ they define a subspace  $\mathbb{R} P^3$
corresponding to the straight line $X \wedge Y$.

On a plane $AP^2$ there are three basis points $E_0, E_1, E_2$
with coordinates
$$
E_0 = (E, 0, 0)^T, \;\; E_1 = (0, E, 0)^T, \;\; E_2 = (0, 0, E)^T,
$$
where $E = \left(
\begin{array}{ll}
1 & 0 \\
0 & 1
\end{array} \right)
$
is the unit matrix, and 0 is the $(2 \times 2)$ 0-matrix.
A point $X \in AP^2$ can be represented in the form
\begin{equation}\label{8}
X = E_0 X^0 + E_1 X^1 + E_2 X^2.
\end{equation}
However, as we noted earlier, the coordinates $X_\alpha$
of this point admit a multiplication from the right
by an element $P \in A$ which is not
a zero divisor.

A point $X$ is in general position with the straight line
$E_\alpha \wedge E_\beta, \, \alpha, \beta = 0, 1, 2,$
if and only if its coordinate $X^\gamma, \, \gamma \neq
\alpha, \beta,$ is not a zero divisor. Let, for instance,
a  point $X$ be in general position with the straight line
$E_1 \wedge E_2$. Then its coordinate $X^0$
is not a zero divisor, and all its coordinates can be
multiplied from the right by $(X^0)^{-1}$. Then expression (8)
of the point $X$ takes the form
\begin{equation}\label{9}
X = E^0 + E_1 \widetilde{X}^1 + E_2 \widetilde{X}^2,
\end{equation}
where $\widetilde{X}^1 = X^1 (X^0)^{-1},
\widetilde{X}^2 = X^2 (X^0)^{-1}$. Now the $(4 \times 2)$-matrix
$(\widetilde{X}^1, \widetilde{X}^2)^T$ is defined uniquely and
is called the \textit{matrix coordinate} of the point $X$
as well as of the straight line $x_0 \wedge x_1$
defined by the point $X$ (see [R 97],  Sect.
{\bf 2.4.1},  and also [R 66], Ch. 3, \S 3).

For the plane $MP^2$, the matrix coordinate has 8 real
components. Hence $\dim \, MP^2 = 8$. Since
$\dim \, MP^2 = \dim \, G (1, 5)$, the plane $MP^2$
can be bijectively mapped onto the Grassmannian $G (1, 5)$.

For the planes $\mathbb{C} P^2,
 \mathbb{C}^1 P^2$, and $\mathbb{C}^0 P^2$,
 the matrix coordinates of points have 4 real components.
 Hence the real dimension of these planes
 is 4,
 $$
\dim \, \mathbb{C} P^2 = \dim \, \mathbb{C}^1 P^2
= \dim \,\mathbb{C}^0 P^2 = 4.
 $$
 Therefore, the family of straight lines $x_0 \wedge x_1$
 in the space $\mathbb{R} P^5$ for each of these planes depends on
4  parameters, i.e., it
forms a \textit{congruence} in the space  $\mathbb{R} P^5$.
We denote these congruences by $K, K^1$, and $K^0$, respectively.

\textbf{4. Moving frames in projective planes over algebras.}
A moving frame in a projective plane $A P^2$ over an algebra $A$
is a triple of points $A_\alpha, \, \alpha = 0, 1, 2,$ that
are mutually not adjacent.
Any point $X \in A P^2$ can be written as
$$
X = A_0 X^0 + A_1 X^1 + A_2 X^2,
$$
where $X^\alpha \in A$ are the coordinates of this point with
respect to the frame $\{A_0, A_1, A_2\}$. The coordinates of a
point $X$ are defined up to a multiplication  from the right
by an element  $P$ of the algebra $A$ which is not a zero divisor.
If a point $X$ is in general position with the straight line $A_1
\wedge A_2$, then its coordinate $X^0$  is not a zero divisor.
Thus, the point $X$ can be written as
$$
X = A^0 + A_1 \widetilde{X}^1 + A_2 \widetilde{X}^2,
$$
where $\widetilde{X}^1 = X^1 (X^0)^{-1},
\widetilde{X}^2 = X^2 (X^0)^{-1}$.
The matrix $(\widetilde{X}^1, \widetilde{X}^2)^T$ is
the matrix coordinate of the point $X$
with respect to the moving frame $\{A_\alpha \}$, and this
coordinate is defined uniquely.

The plane $AP^2$ admits a representation on the Grassmannian
$G (1, 5)$ formed by the straight lines of the space
$\mathbb{R} P^5$. Under this representation, to the vertices
of the  frame $\{A_\alpha \}$ there correspond the straight lines
\begin{equation}\label{10}
A_0 = a_0 \wedge a_1, \;\;
A_1 = a_2 \wedge a_3, \;\;
A_2 = a_4 \wedge a_5
\end{equation}
in $\mathbb{R} P^5$; here $a_i, \, i = 0, \ldots , 5$, are points
of the space $\mathbb{R} P^5$.

The equations of infinitesimal displacement of the
moving frame $\{A_0, A_1, A_2\}$ have the form
\begin{equation}\label{11}
  dA_\alpha = A_\beta \Omega_\alpha^\beta, \;\; \alpha, \beta = 0,
  1, 2,
\end{equation}
where $\Omega_\alpha^\beta$ are  1-forms over the algebra $A$.
In the representation of the algebra $A$ by $(2\times 2)$-matrices,
these forms are expressed as the transposed matrices (1), (2), (3),
and (4). Their entries are not the numbers. The entries are real 1-forms:
\begin{equation}\label{12}
\renewcommand{\arraystretch}{1.5}
\Omega_\alpha^\beta  = \left( \begin{array}{ll}
\omega_{2\alpha}^{2\beta} & \omega_{2\alpha}^{2\beta + 1} \\
\omega_{2\alpha + 1}^{2\beta} & \omega_{2\alpha + 1}^{2\beta + 1}
\end{array} \right).
\renewcommand{\arraystretch}{1}
\end{equation}
Thus, for the plane
$\mathbb{C} P^2$, the entries of the matrix $\Omega_\alpha^\beta$
satisfy the equations
\begin{equation}\label{13}
\omega_{2\alpha}^{2\beta} =  \omega_{2\alpha + 1}^{2\beta + 1}, \;\;
\omega_{2\alpha}^{2\beta + 1} = - \omega_{2\alpha + 1}^{2\beta},
\end{equation}
for the plane $\mathbb{C}^1 P^2$ the equations
\begin{equation}\label{14}
\omega_{2\alpha}^{2\beta} =  \omega_{2\alpha + 1}^{2\beta + 1}, \;\;
\omega_{2\alpha}^{2\beta + 1} = \omega_{2\alpha + 1}^{2\beta},
\end{equation}
and for the plane $\mathbb{C}^0 P^2$ the equations
\begin{equation}\label{15}
\omega_{2\alpha}^{2\beta} =  \omega_{2\alpha + 1}^{2\beta + 1}, \;\;
\omega_{2\alpha + 1}^{2\beta} = 0.
\end{equation}

If now the frame $\{A_\alpha\}$ moves in the plane
$A P^2$, then
the points $a_i \in \mathbb{R} P^5$ also move.
The equations of infinitesimal displacement
of the moving frame  $\{a_i\}$ can be written
in the form
\begin{equation}\label{16}
  da_i = a_j \omega_i^j, \;\; i, j = 0, 1, \ldots, 5,
\end{equation}
where by (10) the forms $\omega_i^j$ coincide
with the corresponding forms (12). The forms $\omega_j^i$
satisfy the structure equations of the projective space
$\mathbb{R} P^5$:
\begin{equation}\label{17}
d \omega_j^i = - \omega_k^i \wedge \omega_j^k,
\end{equation}
where $d$ is the symbol of exterior differential, and
$\wedge$ denotes the exterior multiplication of
the linear differential forms (see for example, [AG 93],
Sec. \textbf{1.3}).

\textbf{5. Focal properties of the congruences
$\boldsymbol{K}, \boldsymbol{K}\boldsymbol{^1},$ and
$\boldsymbol{K}\boldsymbol{^0}$.}
Now we  consider the congruences $K, K^1,$ and
$K^0$ of the space $\mathbb{R} P^5$, representing the planes
$\mathbb{C} P^2, \mathbb{C}^1 P^2,$ and $\mathbb{C}^0 P^2$
in this space, and investigate their
 focal properties.

\begin{theorem}
The projective planes $\mathbb{C} P^2, \mathbb{C}^1 P^2$,
and $\mathbb{C}^0 P^2$ admit a bijective mapping onto the
linear congruences $K, K^1$, and $K^0$ of  the
real space $\mathbb{R} P^5$. These congruences
respectively elliptic, hyperbolic, and parabolic.
\end{theorem}

 {\sf Proof.} We associate to each
 of these congruences a family of projective frames in such
 a way that  the points $a_0$ and $a_1$ are
located on a moving straight line of the congruence.

For the congruence $K$, equations of infinitesimal displacement
of the points $a_0$ and $a_1$ can be written in  the form
\begin{equation}\label{18}
\renewcommand{\arraystretch}{1.5}
\left\{
  \begin{array}{ll}
da_0 = \omega_0^0 a_0  + \omega_0^1 a_1  +\omega_0^2 a_2 +
 \omega_0^3 a_3 +  \omega_0^4 a_4 +  \omega_0^5 a_5, \\
da_1 = - \omega_0^1 a_0  + \omega_0^0 a_1  - \omega_0^3 a_2  +\omega_0^2 a_3
 - \omega_0^5 a_4 +  \omega_0^4 a_5.
  \end{array}
\right.
\renewcommand{\arraystretch}{1}
\end{equation}
By (14), for the congruence  $K^1$, these
two equations  take the form
\begin{equation}\label{19}
\renewcommand{\arraystretch}{1.5}
\left\{
  \begin{array}{ll}
da_0 = \omega_0^0 a_0  + \omega_0^1 a_1  +\omega_0^2 a_2 +
 \omega_0^3 a_3 +  \omega_0^4 a_4 +  \omega_0^5 a_5, \\
da_1 = \omega_0^1 a_0  + \omega_0^0 a_1  +\omega_0^3 a_2  +\omega_0^2 a_3
 + \omega_0^5 a_4 +  \omega_0^4 a_5.
  \end{array}
\right.
\renewcommand{\arraystretch}{1}
\end{equation}
Finally, by (15), for the congruence $K^0$, these
two equations  take the form
\begin{equation}\label{20}
\renewcommand{\arraystretch}{1.5}
\left\{
  \begin{array}{llll}
da_0 = \omega_0^0 a_0  + & \!\!\!\! \omega_0^1 a_1  +\omega_0^2 a_2 &
 \!\!\!\! +  \omega_0^3 a_3 +  \omega_0^4 a_4 & \!\!\!\! +  \omega_0^5 a_5, \\
da_1 =  & \!\!\!\! \omega_0^0 a_1 & \!\!\!\! +\omega_0^2 a_3  & \!\!\!\!  +  \omega_0^4 a_5.
  \end{array}
\right.
\renewcommand{\arraystretch}{1}
\end{equation}

Let $x = a_1 + \lambda a_0$ be an arbitrary point of the straight
line $a_0 \wedge a_1$. This point is a focus of
 this straight line  if for some displacement, its differential
 $dx$ also belongs to this straight line.

 Let us start from the congruence $K^1$, since
 the focal images for this congruence are real and
 look more visual. By (19), for this congruence
 we have
\begin{equation}\label{21}
dx \equiv (\omega_0^3 + \lambda \omega_0^2) a_2
+ (\omega_0^2 + \lambda \omega_0^3) a_3
+ (\omega_0^5 + \lambda \omega_0^4) a_4
+ (\omega_0^4 + \lambda \omega_0^5) a_5 \pmod{a_0 \wedge a_1},
\end{equation}
and as a result, for its focus $x$, the following
equations must be satisfied:
\begin{equation}\label{22}
\renewcommand{\arraystretch}{1.5}
\left\{ \begin{array}{ll}
            \omega_0^2 + \lambda \omega_0^3 = 0, &
            \omega_0^4 + \lambda \omega_0^5 = 0,  \\
 \lambda \omega_0^2 + \omega_0^3 = 0, &
 \lambda \omega_0^4 + \omega_0^5 = 0.
                \end{array} \right.
\renewcommand{\arraystretch}{1}
\end{equation}

The necessary and sufficient
condition of consistency of this system is
$$
\left| \begin{array}{rr}
      1 & \lambda \\
 \lambda & 1
                 \end{array} \right|^2 = 0.
$$
It follows that the values $\lambda = \pm 1$
are double roots of this equation. Thus, each line
$a_0 \wedge a_1$ of  the congruence $K^1$ has two
double foci
$$
f_1 = a_1 + a_0, \;\; f_2 = a_1 - a_0.
$$

Equations (14) imply that the differentials
of the focus $f_1$ are expressed
only in terms of the points $a_0 + a_1,
a_2 + a_3$, and $a_4 +  a_5$. The differentials
of the points $a_2 + a_3$ and $a_4 +  a_5$
are expressed in terms of the same points. As a result, the plane
$$
\pi_1 = (a_0 + a_1) \wedge (a_2 + a_3) \wedge (a_4 +  a_5)
$$
remains fixed when the straight line   $a_0 \wedge a_1$
describes the congruence $K^1$ in the space
$\mathbb{R} P^5$. In a similar way one can prove that
the  focus $f_2$ describes the plane
$$
\pi_2 = (a_0 - a_1) \wedge (a_2 - a_3) \wedge (a_4 -  a_5).
$$

Thus, the congruence $K^1$ is a four-parameter family of straight lines of
the space $\mathbb{R} P^5$ intersecting its two planes $\pi_1$ and
$\pi_2$ that are in general position. Hence $K^1$
 is a  \textit{hyperbolic} line congruence.

In a similar way, we can prove that each straight line
$a_0 \wedge a_1$ of the congruence $K$ bears two double complex
conjugate foci
$$
f_1 = a_1 + i a_0, \;\; f_2 = a_1 - i a_0,
$$
and these foci describe two
 complex conjugate two-dimensional planes
 $\pi_1$ and $\pi_2, \pi_2 = \overline{\pi}_1$.
Hence $K$  is an  \textit{elliptic} line congruence
in the space $\mathbb{R} P^5$. The straight lines
of $K$ do not have singular points in  $\mathbb{R} P^5$.

Consider finally the  congruence $K^0$ in
the space $\mathbb{R} P^5$. We will look for
the foci of its straight lines in the same form
$$
 x = a_1 + \lambda a_0.
$$
Differentiating this expression by means of (20), we find that
$$
dx \equiv \lambda \omega_0^2 a_2
+ (\lambda \omega_0^3 + \omega_0^2) a_3
+ \lambda \omega_0^4 a_4
+ (\lambda \omega_0^5 + \omega_0^4) a_5 \pmod{a_0 \wedge a_1}.
$$
Thus the focus $x$ must  satisfy the following
equations:
\begin{equation}\label{23}
\renewcommand{\arraystretch}{1.5}
\left\{ \begin{array}{llllll}
            \lambda \omega_0^2  &\!\!\!\! &\!\!\!\! =  0, &
             \lambda \omega_0^4  &\!\!\!\!  &\!\!\!\! = 0,  \\
  \omega_0^2  &\!\!\!\!+ \lambda \omega_0^3  &\!\!\!\!= 0, &
  \omega_0^4  &\!\!\!\!+ \lambda \omega_0^5  &\!\!\!\!= 0.
                \end{array} \right.
\renewcommand{\arraystretch}{1}
\end{equation}
This system is consistent if and only if
$$
\left| \begin{array}{rr}
      \lambda & 0 \\
 1 & \lambda  \end{array} \right|^2 = 0.
$$
It follows that the value $\lambda = 0$
is a quadruple root of this equation. Thus, each line
$a_0 \wedge a_1$ of  the congruence $K^0$ has a real quadruple
singular point $f = a_1$. Applying equations (15), it is easy to
prove that when the straight line
$a_0 \wedge a_1$ describes the congruence $K^0$, this
focus describes the plane $\pi = a_1 \wedge a_3 \wedge a_5$.
Hence $K^0$  is a  \textit{parabolic} line congruence.
\rule{3mm}{3mm}

\textbf{6. Smooth lines in projective planes.} On a
projective plane $A P^2$, where $A$ is one of the algebras
$\mathbb{C}, \mathbb{C}^1$, and $\mathbb{C}^0$, consider a
smooth point submanifold $\gamma$ of real dimension 2.
Such a submanifold is called an \textit{$A$-smooth line}
if at any of its points $X$, it is tangent to a straight line
$U$ passing through $X$.

With an $A$-smooth line $\gamma$, associate a family
of projective frames $\{A_0 A_1 A_2\}$ in such a way
that $A_0 = X$ and $A_1$ lies on the tangent $U$ to
$\gamma$ at $X$. Then on the line $\gamma$, the first
equation of (11) takes the form
\begin{equation}\label{24}
d A_0 = A_0 \Omega_0^0 + A_1 \Omega_0^1.
\end{equation}
It follows that $A$-smooth lines on a plane $A P^2$
are defined by the equation
\begin{equation}\label{25}
  \Omega^2_0 =  0.
\end{equation}
The 1-form $\Omega^1_0$ in equation (24) defines
a displacement of the point $A_0$ along the
curve $\gamma$. So, this form is a basis form on $\gamma$.

By equations (12) we have
$$
\renewcommand{\arraystretch}{1.5}
\Omega^1_0 =  \left( \begin{array}{ll}
\omega_{0}^{2} & \omega_{0}^{3}\\
\omega_1^2  & \omega_{1}^{3},
\end{array} \right), \;\;
\Omega^2_0 =  \left( \begin{array}{ll}
\omega_{0}^{4} & \omega_{0}^{5} \\
\omega_1^4 & \omega_{1}^{5}
\end{array} \right),
\renewcommand{\arraystretch}{1}
$$
where $\omega_i^j$ are  real 1-forms.
For the algebras
$\mathbb{C}, \mathbb{C}^1$, and $\mathbb{C}^0$,
they are related respectively by equations (13), (14), and (15).
As a result, on the line $\gamma \subset A P^2$, the following
differential equations will be satisfied:
\begin{equation}\label{26}
  \omega^4_0 =  0, \;\; \omega^5_0 =  0.
\end{equation}
These equations are equivalent to equations (25).

Since $\Omega_0^1$ is a basis form on the plane $A P^2$,
the real forms $\omega_0^2$ and  $\omega_0^3$ are
linearly independent.  The families of straight lines in the space
$\mathbb{R} P^5$ corresponding to these lines depend on 2
parameters and form a three-dimensional ruled submanifold.
Denote this submanifold by $S$. These submanifolds
coincide with the congruences $K, K^1$, and $K^0$ if
$\gamma \subset \mathbb{C} P^2, \gamma \subset \mathbb{C}^1 P^2$,
and $\gamma \subset \mathbb{C}^0 P^2$, respectively.

\begin{theorem}
The tangent subspace $T_x (S)$ to the ruled submanifold $S$
corresponding in the space $\mathbb{R} P^5$ to a smooth line
in the planes  $\mathbb{C} P^2, \mathbb{C}^1 P^2$,
and $\mathbb{C}^0 P^2$  is fixed at all points of its
rectilinear generator $L$, and
the submanifold $S$ is tangentially degenerate of rank $r = 2$.
\end{theorem}

{\sf Proof.}
Consider a rectilinear generator $L = a_0 \wedge a_1$
of the submanifold $S$. By (26), the differentials of the points
$a_0$ and $a_1$ are written in the form
\begin{equation}\label{27}
\renewcommand{\arraystretch}{1.5}
\left\{
  \begin{array}{ll}
da_0 = \omega_0^0 a_0  + \omega_0^1 a_1  +\omega_0^2 a_2 +
 \omega_0^3 a_3, \\
da_1 = \omega_0^1 a_0  + \omega_1^1 a_1  +\omega_1^2 a_2  +\omega_1^3 a_3.
  \end{array}
\right.
\renewcommand{\arraystretch}{1}
\end{equation}
It follows that at any point
$x \in a_0 \wedge a_1$, the tangent subspace $T_x (X)$
belongs to a three-dimensional subspace $P^3 \subset \mathbb{R}
P^5$ defined by the points $a_0, a_1, a_2$, and $a_3$.
Thus the subspace $T_x (X)$ remains fixed along
the rectilinear generator $L = a_0 \wedge a_1$, and
the submanifold $S$ is tangentially degenerate of rank $r = 2$.
\rule{3mm}{3mm}

\textbf{7. Singular points of submanifolds corresponding to smooth
lines in the projective spaces over two-dimensional algebras.}
We will prove  the following theorem.

\begin{theorem}
To smooth lines in the projective planes
$\mathbb{C} P^2, \mathbb{C}^1 P^2$, and $\mathbb{C}^0 P^2$
over the algebras of complex, double, and dual numbers,
there correspond three-dimensional tangentially
degenerate submanifolds of rank $r = 2$
in the space $\mathbb{R} P^5$. For the algebra
$\mathbb{C}$ such a submanifold does not have
real singular points, for the algebra
$\mathbb{C}^1$, such a submanifold is the
join formed by the straight lines connecting
the points of two planes curves that are in general
position, and  for the algebra
$\mathbb{C}^0$, such a submanifold is
a subfamily of the family of straight lines
intersecting a plane curve. In all these cases,
a submanifold $S$  depends on two functions of one variable.
\end{theorem}

A rectilinear generator $L = a_0 \wedge a_1$
of a submanifold $S$ of rank 2 bears two foci. Let us find these
foci for the submanifolds $S$ corresponding to the
lines $\gamma$ in the planes
$\mathbb{C} P^2, \mathbb{C}^1 P^2$, and $\mathbb{C}^0 P^2$.
We assume that  these foci have the form $x = a_1 + \lambda a_0$.

If a line $\gamma \subset \mathbb{C}^1 P^2$, then
equations (14) and (27) are satisfied. They imply that
$$
dx \equiv (\omega_0^3 + \lambda \omega_0^2) a_2
+ (\omega_0^2 + \lambda \omega_0^3) a_3
\pmod{a_1 \wedge a_2},
$$
and for the focus $x$, we have
$$
\omega_0^3 + \lambda \omega_0^2 = 0, \;\;
\omega_0^2 + \lambda \omega_0^3 = 0.
$$
This system is consistent if and only if
$$
\left| \begin{array}{rr}
   1 & \lambda  \\
  \lambda & 1 \end{array} \right| = 0,
$$
i.e., if $\lambda = \pm 1$. Thus, the foci of
the straight line $a_0 \wedge a_1$ are
the points $a_1 + a_0$ and $a_1 - a_0$. These
points belong to the focal planes
$\pi_1$ and $\pi_2$ of the congruence
$K^1$ and describe lines $\gamma_1$ and $\gamma_2$.
Such manifolds $S$ are called \textit{joins}.
Since each of the lines $\gamma_1$ and $\gamma_2$
on the  planes $\pi_1$ and $\pi_2$ is defined
by means of one function of one variable, a submanifold $S$
depends on two functions of one variable. The same result
could be obtained by applying the Cartan test
(see [BCGGG 90]) to the system of equations (14) and (26).

If a line $\gamma \subset \mathbb{C} P^2$, then
we can prove that a rectilinear generator $L = a_0 \wedge a_1$
of the  ruled submanifold $S$ corresponding to $\gamma$ bears
two complex conjugate foci belonging to
complex conjugate focal planes $\pi_1$ and $\pi_2 = \overline{\pi}_1$
of the congruence $K$. Hence \textit{in the real space
$\mathbb{R} P^5$, the submanifold $S$ does not have singular
points.}

In the complex plane $\pi_1$, the focus $f_1$ can describe
an arbitrary differentiable line. But such a line
is defined by means of two functions of one real variable.
Therefore, in this case the submanifold $S$ also
depends on  two functions of one real variable.

Finally, consider a submanifold $S \subset \mathbb{R} P^5$
corresponding to a line $\gamma \subset \mathbb{C}^0 P^2$.
Such a submanifold is defined in $ \mathbb{R} P^5$
by differential equations (15) and (26). Using
the same method as above, we can prove that
a rectilinear generator $L = a_0 \wedge a_1$
of the  ruled submanifold $S$ corresponding to $\gamma$ bears
a double real  focus $f = a_1$ belonging to
the focal plane $\pi$ of the congruence $K^0$
and describing in this plane an arbitrary line.

We prove that in this case
a submanifold $S$ is also defined by two functions
of one variable. Now we will apply the Cartan test.

Taking exterior derivatives of equations (26) and applying
equations (15),  we obtain the following exterior
quadratic equations:
\begin{equation}\label{28}
  \omega^2_0 \wedge \omega^4_2 =  0, \;\;
   \omega^2_0 \wedge \omega^5_2 + \omega^3_0 \wedge \omega^4_2 =  0.
\end{equation}
It follows from (28) that
\begin{equation}\label{29}
   \omega^4_2 =  a \omega_0^2, \;\;
   \omega^5_2 =  b \omega_0^2 + a \omega_0^3.
\end{equation}

 We apply the Cartan test to
the system of equations (26), (28), and (29).
In addition to the basis forms $\omega_0^2$
and $\omega_0^3$, equations (28) contain
two more forms $\omega_2^4$ and $\omega_2^5$. Thus,
we have $q = 2$. The number of independent equations
in (28) is also 2, i.e., $s_1 = 2$. As a result,
$s_2 = q - s_1 = 0$, and the Cartan number
$$
Q = s_1 + 2 s_2 = 2.
$$
Equations (29) show that the number $N$ of parameters
on which the general two-dimensional integral element depends
is also 2, $N = 2$. Since $Q = N$, the system of equations (26)
is involution, and its solution  depends on two functions
of one variable. \rule{3mm}{3mm}

Note also that in the space $\mathbb{R} P^5$,
to the straight lines in the planes $\mathbb{C} P^2,
\mathbb{C}^1 P^2$, and $\mathbb{C}^0 P^2$ there correspond
 two-parameter linear congruences of straight lines
 located in its three-dimensional subspaces.
For the plane $\mathbb{C} P^2$, these congruences are elliptic,
for the plane $\mathbb{C}^1 P^2$, they are hyperbolic,
and for the plane $\mathbb{C}^0 P^2$, they are parabolic.

\vspace*{10mm}

\noindent {\em Authors' addresses}:\\

\noindent
\begin{tabular}{ll}
M.~A. Akivis &V.~V. Goldberg\\
Department of Mathematics &Department of Mathematical Sciences\\
Jerusalem College of Technology---Mahon Lev &  New
Jersey Institute of Technology \\
Havaad Haleumi St., P. O. B. 16031 & University Heights \\
 Jerusalem 91160, Israel &  Newark, N.J. 07102, U.S.A. \\
 & \\
 E-mail address: akivis@avoda.jct.ac.il & E-mail address:
 vlgold@m.njit.edu
 \end{tabular}
\end{document}